\def\timestamp{%
Time-stamp: <1st-ctble.tex: Friday 30-05-2008 at 14:21:32 (cest)>}
\def\stripname Time-stamp: <#1 #2>{#2}
\edef\filedate{\expandafter\stripname\timestamp}

\documentclass[a4paper]{amsart}

\usepackage{amssymb,amsrefs}

\newcommand{\Ex}{\operatorname{Ex}}
\newcommand{\Fn}{\operatorname{Fn}}
\newcommand{\orpr}[2]{\langle#1,#2\rangle}
\newcommand{\preim}{^\gets}

\DeclareMathSymbol\forces     \mathrel{AMSa}{"0D}
\DeclareMathSymbol\restr2{AMSa}{"16}
\DeclareMathSymbol\HH0{AMSb}{`H} 
\DeclareMathSymbol\N0{AMSb}{`N} 
\DeclareMathSymbol\Q0{AMSb}{`Q} 
\DeclareMathSymbol\R0{AMSb}{`R} 

\newcommand\betaH{\beta\HH}

\newcommand\Hstar{\HH^*}
\newcommand\Nstar{\N^*}
\newcommand\CE{C_E}

\let\phi\varphi

\newcommand{\CH}{\mathsf{CH}}

\newcommand{\Pow}{\mathcal{P}}

\theoremstyle{remark}

\newtheorem*{remark}{Remark}

\begin{document}

\title{A first-countable non-remainder of $\HH$}

\author[Alan Dow]{Alan Dow\dag}
\address{
Department of Mathematics\\
UNC-Charlotte\\
9201 University City Blvd. \\
Charlotte, NC 28223-0001}
\email{adow@uncc.edu}
\urladdr{http://www.math.uncc.edu/\~{}adow}
\thanks{\dag Supported by NSF grant DMS-0554896}

\author{Klaas Pieter Hart}
\address{Faculty of Electrical Engineering, Mathematics and Computer Science\\
         TU Delft\\
         Postbus 5031\\
         2600~GA {} Delft\\
         the Netherlands}
\email{k.p.hart@tudelft.nl}
\urladdr{http://fa.its.tudelft.nl/\~{}hart}

\date{\filedate}

\begin{abstract}
We give a (consistent) example of a first-countable continuum that 
is not a remainder of the real line.
\end{abstract}

\subjclass[2000]{Primary: 54F15. 
                 Secondary: 03E50, 03E65, 54A35, 54D40, 54G20}
\keywords{first-countable continuum, continuous image, 
          $\Hstar$, Cohen reals}

\maketitle

\section*{Introduction}

The purpose of this note is to confirm a suspicion raised 
in~\cite{dowhart-sep}*{Question~4.2}:
we show that Bell's example, from~\cite{MR1058795}, of a first-countable
compact space that is not an $\Nstar$-image can be adapted to produce
admits a connected variation that is neither an $\Nstar$-image nor an 
$\Hstar$-image.
The interest in this variation stems from the authors' 
version of Parovi\v{c}enko's theorem from~\cite{Parovicenko63}.
That theorem states that every compact Hausdorff space of weight~$\aleph_1$
or less is an $\Nstar$-image; the Continuum Hypothesis then implies that
the $\Nstar$-images are exactly the compact Hausdorff spaces of 
weight~$2^{\aleph_0}$ or less.
We proved in~\cite{MR1707489} a parallel result for $\Hstar$ and continua 
(connected compact Hausdorff spaces).
Since, by Arkhangel\cprime ski\u{\i}'s theorem~\cite{MR0251695},
first-countable compact spaces have weight at most~$2^{\aleph_0}$ 
it follows that under $\CH$ first-countable compacta/continua are 
$\Nstar$-images/$\Hstar$-images respectively.

\subsection*{Bell's graph}

A major ingredient in our construction is Bell's graph, constructed
in~\cite{MR677860}.
It is a graph on the ordinal~$\omega_2$, represented by a symmetric
subset~$E$ of~$\omega_2^2$.
The crucial property of this graph is that there is \emph{no} map 
$\varphi:\omega_2\to\Pow(\N)$ that represents this graph in the sense that
$\orpr\alpha\beta\in E$ if and only if $\varphi(\alpha)\cap\varphi(\beta)$ is
infinite.

Bell's graph exists in any forcing extension in which $\aleph_2$~Cohen 
reals are added; for the reader's convenience we shall describe the
construction of~$E$ and adapt Bell's proof so that it applies to continuous 
maps defined on~$\Hstar$.

\section*{A first-countable continuum}

Our starting point is a connected version of the Alexandroff double of the
unit interval.
We topologize the unit square as follows.
\begin{enumerate}
\item a local base at points of the form $\orpr x0$ consists of the sets
      $$
       U(x,0,n)=(x-2^{-n},x+2^{-n})\times[0,1] \setminus
                \{x\}\times[2^{-n},1]
      $$
\item a local base at points of the form $\orpr xy$, with $y>0$ consists
      of the sets
      $$
       U(x,y,n) = \{x\}\times(y-2^{-n},y+2^{-n})
      $$
\end{enumerate}
We call the resulting space the connected comb and denote it by~$C$.
It is straightforward to verify that $C$ is compact, Hausdorff and connected;
it is first-countable by definition.

For each $x\in[0,1]$ and positive $a$ we define to be the following
cross-shaped closed subset of~$C^2$:
$$
D_{x,a}= 
\bigl(\{x\}\times[a,1]\times C\bigr)\cup\bigl(C\times\{x\}\times[a,1]\bigr)
$$
We note the following two properties of the sets $D_{x,a}$
\begin{enumerate}
\item if $a<b$ then $D_{x,b}$ is in the interior of~$D_{x,a}$, and
\item if $x\neq y$ then $D_{x,a}\cap D_{y,a}$ is the union of two squares:
       $\{x\}\times[a,1]\times\{y\}\times[a,1]$ and 
       $\{y\}\times[a,1]\times\{x\}\times[a,1]$
\end{enumerate}

Now take any $\aleph_2$-sized subset of~$[0,1]$ and index it (faithfully)
as $\{x_\alpha:\alpha<\omega_2\}$.
We use this indexing to identify $E$ with the subset 
$\{\orpr{x_\alpha}{x_\beta}:\orpr\alpha\beta\in E\}$ of the unit square.
Next we remove from $C^2$ the following open set:
$$
\bigcup_{\orpr xy\notin E}
\Bigl(\bigl(\{x\}\times(0,1]\times\{y\}\times(0,1]\bigr) \cup
\bigl(\{y\}\times(0,1]\times\{x\}\times(0,1]\bigr)\Bigr)
$$
The resulting compact space we denote by $\CE$.
Observe that the intersections $D_{x_\alpha,a}\cap \CE$ represent~$E$ in the 
sense that $D_{x_\alpha,a}\cap D_{x_\beta,a}\cap \CE$ is nonempty if and only
if~$\orpr\alpha\beta\in E$.
We write $D^E_{x,a}=D_{x,a}\cap \CE$.

\smallskip

We show that $\CE$ is (arcwise) connected.  

To begin: the square~$S$ of the base line of~$C$ is a subset of~$\CE$ and 
homeomorphic to the unit square so that it is (arcwise) connected.

Let $\langle x,a,y,b\rangle$ be a point of~$\CE$ not in~$S$.
If, say, $a=0$ then $\{\orpr x0\}\times(\{y\}\times[0,b])$ is an arc 
in~$\CE$ that connects $\langle x,0,y,b\rangle$ to the 
point~$\langle x,0,y,0\rangle$ in~$S$.
If $a,b>0$ then $\orpr xy\in E$ and the whole square
$\{x\}\times[0,1]\times\{y\}\times[0,1]$ is in~$\CE$ and it provides
us with an arc in~$\CE$ from $\langle x,a,y,b\rangle$ 
to $\langle x,0,y,0\rangle$.

\smallskip

We find that $\CE$ is a first-countable continuum.
It remains to show that it is not an $\Hstar$-image.

Assume $h:\Hstar\to \CE$ is a continuous surjection and consider, 
for each~$\alpha$, the sets $D^E_{x_\alpha,\frac34}$ and $D^E_{x_\alpha,\frac12}$.

Using standard properties of $\betaH$, see~\cite{Hart}*{Proposition~3.2},
we find for each~$\alpha$ a sequence 
$\bigl<(a_{\alpha,n},b_{\alpha,n}):n\in\N\bigr>$
of open intervals with rational endpoints, 
and with $b_{\alpha,n}<a_{\alpha,n+1}$ for all~$n$, 
such that 
$h\preim[D^E_{x_\alpha,\frac34}]\subseteq 
    \Ex O_\alpha\cap\Hstar \subseteq h\preim[D^E_{x_\alpha,\frac12}]$,
where $O_\alpha=\bigcup_n(a_{\alpha,n},b_{\alpha,n})$.
Because the intersections of the sets~$D^E_{x_\alpha,a}$ represent~$E$
the intersections of the $O_\alpha$ will do this as well:
the conditions `$O_\alpha\cap O_\beta$ is unbounded' 
and `$\orpr\alpha\beta\in E$' are equivalent.

In the next subsection we show that for (many) $\orpr\alpha\beta$
this equivalence does not hold and that therefore $\CE$ is not a continuous
image of~$\Hstar$.

Note also that our continuum is not an $\Nstar$-image either: 
if $g:\Nstar\to \CE$ were continuous and onto we could use clopen subsets
of~$\Nstar$ and their representing infinite subsets of~$\N$ to contradict
the unrepresentability property of~$E$.

\subsection*{Destroying the equivalence}

We follow the argument from~\cite{MR677860} and we rely on Kunen's book 
\cite{MR597342}*{Chapter~VII} for basic facts on forcing.
We let $L=\{\orpr\alpha\beta\in\omega_2^2:\alpha\le\beta\}$ 
and we force with the partial order $\Fn(L.2)$ of finite partial functions 
with domain in~$L$ and range in~$\{0,1\}$.
If $G$~is a generic filter on~$\Fn(L,2)$ then we let 
$E=\{\orpr\alpha\beta: \bigcup G(\alpha,\beta)=1$ or 
$\bigcup G(\beta,\alpha)=1\}$.

To show that $E$ is as required we take a nice name~$\dot F$ for a function 
from~$\omega_2$ to~$(\Q^2)^\omega$ that represents a choice of open 
sets~$\alpha\mapsto O_\alpha$ as in above in that  
$F(\alpha)=\bigl<\orpr{a_{\alpha,n}}{b_{\alpha,n}}:n\in\omega\bigr>$
for all~$\alpha$.
As a nice name $\dot F$ is a subset of 
$\omega_2\times\omega\times\Q^2\times\Fn(L,2)$, where for each point
$\langle\alpha,n,a,b\rangle$ the set 
$\{p:\langle\alpha,n,a,b,p\rangle\in\dot F\}$ is a maximal antichain
in the set of conditions that forces the $n$th term of $\dot F(\alpha)$
to be $\orpr ab$.

For each $\alpha$ we let $I_\alpha$ be the set of ordinals that occur in the
domains of the conditions that appear as a fifth coordinate in the elements 
of~$\dot F$ with first coordinate~$\alpha$.
The sets $I_\alpha$ are countable, by the ccc of $\Fn(L,2)$.
We may therefore apply the Free-Set Lemma,
see \cite{MR795592}*{Corollary~44.2},
and find a subset~$A$ of~$\omega_2$ of cardinality~$\aleph_2$ such that
$\alpha\notin I_\beta$ and $\beta\notin I_\alpha$ whenever 
$\alpha,\beta \in A$ and $\alpha\neq\beta$.

Let $p\in\Fn(L,2)$ be arbitrary and take $\alpha$ and~$\beta$ in~$A$ with 
$\alpha<\beta$ and such that $\alpha>\eta$ whenever $\eta$ occurs in~$p$.
Consider the condition $q=p\cup\bigl\{\langle\alpha,\beta,1\rangle\bigr\}$. 
If $q$~forces $O_\alpha\cap O_\beta$ to be bounded in~$[0,\infty)$ then we 
are done: $q$~forces that the equivalence fails at $\orpr\alpha\beta$.

If $q$~does not force the intersection to be bounded we can extend~$q$
to a condition~$r$ that forces $O_\alpha\cap O_\beta$ to be unbounded.
We define an automorphism $h$ of~$\Fn(L,2)$ by changing the value of the
conditions only at~$\orpr\alpha\beta$: from~$0$ to~$1$ and vice versa.
The condition~$p$ as well as the names~$\dot x_\alpha$ and $\dot x_\beta$
are invariant under~$h$.
It follows that $h(r)$ extends~$p$ and
$$
h(r)\forces \hbox{$\bigcup\dot G(\alpha,\beta)=0$ and 
   $O_\alpha\cap O_\beta$ is unbounded}
$$
so again the equivalence is forced to fail at~$\orpr\alpha\beta$.

\begin{remark}
The argument above goes through almost verbatim to show that Bell's graph
can also be obtained adding $\aleph_2$ random reals.
When forcing with the random real algebra one needs only consider condtions
that belong to the $\sigma$-algebra generated by the clopen sets of the
product $\{0,1\}^L$; these all have countable supports so that, 
again by the ccc, one can define the sets $I_\alpha$ as before.
The rest of the argument remains virtually unchanged.
\end{remark}

\end{document}